\def\IH{{\Bbb H}} 
 
\def\IR{{\Bbb R}} 
 
\def\IS{{\Bbb S}}

\def\IC{\Bbb C} 
\def\ID{{\Bbb D}}
 
\def\tOmega{\tilde{\Omega}} 
 
\def\zbar{{\overline{z}}} 
 
\def\wbar{{\overline{w}}} 
\def\mubar{{\overline{\mu}}} 

\documentclass[12pt]{article} 

\usepackage[psamsfonts]{amssymb} 
\usepackage{amsfonts,amsmath}

\usepackage{epsfig,multicol}

\newtheorem{theorem}{Theorem}
\newtheorem{lemma}{Lemma}
\newtheorem{corollary}{Corollary}

\title{The tension equation with holomorphic coefficients,  harmonic mappings and rigidity} 
\author{Gaven J. Martin \thanks{g.j.martin@massey.ac.nz \newline
New Zealand Institute for Advanced Study,  Massey University, NZ and \newline Magdalen College,  Oxford University \newline Research supported in parts by 
grants from the  N.Z. Marsden Fund. }} 
\date{} 
\begin{document}

\maketitle 


\begin{abstract} \noindent The tension equation for a mapping $f:\IC\to\IC$ is the nonlinear second order equation
\[ \Delta f +\varphi(f) f_z f_\zbar = 0\]
Solutions are ``harmonic'' mappings.  Here we give a complete description of the solution space of mappings of degree $1$ to this equation when $\varphi$ is entire.  Each solution is a quasiconformal surjection and when the set of normalised solutions is endowed with the Teichm\"uller metric,  the solution space is isometric to the hyperbolic plane.  More generally,  for  harmonic mappings $f:\Omega \to (\tilde{\Omega},\rho)$ between domains in  $\IC$, with $\rho(w)|dw|$ defining a flat metric we stablish a very strong maximum principle for the distortion - up to multiplicative  factor $e^{iv}$,  $v$ real and harmonic,  the Beltrami coefficient of $f^{-1}$ is quasiregular - and thus open and discrete when nonconstant.  This follows from the remarkable fact that the Beltrami coefficient of the inverse of a harmonic mapping itself satisfies a nonlinear homogeneous Beltrami  equation.   
\end{abstract}

\section{Introduction} Let $\Omega$ and $\tOmega$ be domains in the complex plane $\IC$.   A smooth real valued function $\rho:\tOmega\to\IR_+$ defines a conformal  metric $\rho(w) |dw|$ on $\Omega$.  The Gaussian curvature of this metric is given pointwise by ${\cal K}  = -\rho^{-2} \, \Delta \log \rho$.  The metric is flat if ${\cal K}\equiv 0$,  equivalently if $\log \rho$ is a harmonic real valued function. 

A mapping $f:\Omega \to (\tOmega,\rho)$ is   harmonic  if it is a solution to the nonlinear second order equation, called the {\em tension equation}, 
\begin{equation}\label{te}
f_{z\zbar}(z) + (\log \rho)_w(f)\; f_z(z)\; f_\zbar(z) = 0
\end{equation} 
Note that this criterion does not depend on the conformal structure of the domain of the mapping.
We always assume that our mappings are orientation preserving.   In \cite{Mdiff} we proved that the Jacobian determinant of a harmonic homeomorphism does not vanish, thus harmonic homeomorphisms are diffeomorphisms,  generalising earlier results of Lewy \cite{Lewy} and of Schoen-Yau, \cite{SY}.  Much of the basic material concerning harmonic maps can be found in Jost's book \cite{Jost}.

\medskip

One of the strongest metrics one can put on a space of  mappings is the Teichm\"uller metric;  if $f,g:\Omega\to\tOmega$ define
\[ d_{{\bf T}}(f,g) = \log \sup_{z\in g(\Omega)} K(z,f\circ g^{-1}) = \log \sup_{z\in f(\Omega)} K(z,g\circ f^{-1}) \]
Here $K(z,h)$ is the distortion of a mapping, 
\[ K(\cdot,h) = \frac{|h_z|+|h_\zbar|}{|h_z|-|h_\zbar|} \]
Convergence is this topology implies local uniform convergence, or if $\Omega = \IC$ uniform convergence in the spherical metric.  One of our more interesting results here is the following.

\begin{theorem}\label{isothm}  Let ${\cal F}^{\, 0}=\{f:\IC\to(\IC,\rho)\}$ be the family of entire harmonic homeomorphisms with respect to a flat metric $\rho(w)|dw|$ and normalised by $f(0)=0$,  $f(1)=1$.  Then,  with the Teichm\"uller metric,  $({\cal F}^{\, 0},d_{{\bf T}})$ is canonically isometric to the hyperbolic plane and every $f\in {\cal F}^{\, 0}$ is quasiconformal.
\end{theorem}
The normalisation here is to remove the obvious action of the similarity group ${\cal S}=\{\alpha  z+\beta:\alpha,\beta\in \IC\}$.  The reader will see that Theorem \ref{isothm} really describes the space of normalised homeomorphic (more generally,  degree $1$) solutions to the equation
\begin{equation}\label{te2}
f_{z\zbar} + \Psi(f) \; f_z\; f_\zbar = 0
\end{equation} 
for entire holomorphic $\Psi$.  This is because $\rho$ flat implies $\log \rho$ harmonic and $(\log \rho)_w$ holomoprhic.  Conversely,  if $\Psi$ is entire it admits an antiderivative,  say $U+iV$ with $U_z=\Psi$, and with $\rho=e^U$ we have $\Psi = (\log e^U)_z$.
\begin{corollary}   Let ${\cal F}^{\, 0}=\{f:\IC\to(\IC,\rho)\}$ be the family of degree 1 solutions to the tension equation (\ref{te2}) with $\Psi$ entire and normalised by $f(0)=0$,  $f(1)=1$.  Then,   $({\cal F}^{\, 0},d_{{\bf T}})$ is canonically isometric to the hyperbolic plane.
\end{corollary}

\medskip

\noindent{\bf Example:} In the most elementary case,  $\rho(z)\equiv 1$ or $\Psi\equiv 0$, Berstein's theorem on miminal graphs (or more direct argument, of course) has the consequence that an entire harmonic mapping of the plane $\IC$ is linear, $z\mapsto az+b\zbar+c$,  $|a|> |b|$.  The normalisation puts $c=0$ and $a+b=1$ and $\mu_f = b/a = 1/a-1 \in \ID$ and 
\[ K(f)=\frac{1+|\alpha|}{1-|\alpha|} = \frac{|a|+|1-a|}{|a|-|1-a|} \]
Next if $f=az+(1-a)\zbar$ and $g   =  bz+(1-b)\zbar$,  then $g^{-1} = (-\bar b w+(1-b)\bar w)/(|1-b|^2-|b|^2)$.  Hence 
\[f\circ g^{-1}(w)  =   \frac{(1-a-\bar b) w+ (a-b) \bar w }{|1-b|^2-|b|^2} \]
so that 
\[ K(f\circ g^{-1}) =  \frac{|1-a-\bar b| + |a-b|}{|1-a-\bar b|-|a-b|} \]
  Next,  with $d_\ID$ denoting the hyperbolic metric of the disk, 
\[
d_\ID(\mu_f,\mu_g)=d_{\ID}(\frac{1-b}{b},\frac{1-a}{a}) =    d_\ID\Big(0, \frac{|a-b|}{|1-a-\bar b|}\Big) = \log  K(f\circ g^{-1}) \\
\]
exhibiting the isometry;  harmonic quasiconformal $f \leftrightarrow \mu_f\in \ID$.  For other metrics the Beltrami coefficient will not be constant,  but it will have constant modulus - reinforcing the connection with Teichm\"uller mappings.

\section{Harmonic homeomorphisms}

We now give the basic definitions and describe some results we will need.
 \subsection{The Hopf differential}

The Hopf differential of a harmonic mapping $f$ is
\begin{equation}
\Phi_f(z) = \rho^2(f(z)) f_z(z) \overline{f_\zbar(z)}.
\end{equation}
It is well known that the  Hopf differential of a harmonic mapping is holomorphic.  There are  strong refinements of the converse statement,  see \cite{Iwaniecetc}.   

\subsection{Quasiconformality}

For a harmonic homeomorphism $f:\Omega\to\tOmega$  we define the  the complex dilatation (or Beltrami coefficient) of $f$ as
\begin{equation}
\mu^f(z) = \frac{f_\zbar(z)}{f_z(z)}
\end{equation}
Since $f_z\neq 0$,   $\mu^f\in C^1(\Omega)$.   The distortion of $f$ is then
\begin{equation}
K(z,f)=  \frac{1+|\mu^f(z)|}{1-|\mu^f(z)|}, 
\end{equation}
and $f$ is {\em quasiconformal} if   
$ \|K(z,f)\|_{L^\infty(\ID)} < \infty$.

\medskip

We refer to \cite{AIM} for the basic properties of quasiconformal mappings.   
 \section{Subharmonicity of $\log |\mu_f|$}
 
Let $\Phi_f$ be the holomorphic Hopf differential of the harmonic map $f:\Omega \to (\tOmega,\rho)$.  Fix a point $z_0 \in \Omega\setminus \{z:f_\zbar(z) = 0 \}$.  We choose a well defined branch of the argument in a   domain containing  $f_z(z_0)$ and $f_\zbar(z_0)$ so that we can work with well defined functions $\log f_z = \log |f_z|+i \arg(f_z)$ and also $\log f_\zbar= \log |f_\zbar|+i \arg(f_\zbar)$.  Notice that as $\phi_f=\rho^2(f) f_z \overline{f_\zbar}$ is holomorphic,  $\arg(f_z \overline{f_\zbar})$ is harmonic.  This is basically why the imaginary terms disappear in the following calculation.    For a harmonic mapping $f:\Omega \to (\tOmega,\rho)$,  with $\rho$ flat,  we have
   \begin{eqnarray*}
\lefteqn{ \frac{1}{4}\big(       \Delta   \log   {  f_\zbar} - \Delta   \log   {f_z } \big) =    \Big[\frac{ {  f_{z\zbar}}}{ {  f_{\zbar}} } \Big]_\zbar - \Big[ \frac{{f_{z\zbar} }}{{f_{z} }} \Big]_z 
= \Big[  {{(\log \rho)_w(f)   f_\zbar }} \Big]_z - \Big[  { (\log \rho)_w(f) f_z  } \Big]_\zbar}\\
   & = &  
  \Big[  {(\log \rho)_w(f)\Big]_z   {f_\zbar }} + {(\log \rho)_w(f)  {f_{z\zbar} }}  -  \Big[   (\log \rho)_w (f) \Big]_\zbar {  f_z  }  -(\log \rho)_w (f)   {  f_{z\zbar}  }  \\
& = &    \Big[  {(\log \rho)_w(f)\Big]_z    {f_\zbar }}  -  \Big[   (\log \rho)_w (f) \Big]_\zbar {  f_z  }   = \Big[  {(\log \rho)_{ww}(f) f_z + (\log \rho)_{w\wbar}(f)\overline{f_\zbar} \Big]    {f_\zbar }}  \\ && -  \Big[   (\log \rho)_{ww} (f) f_\zbar + (\log \rho)_{w\wbar} (f) \overline{f_z} \Big]{  f_z  }  = 0  
 \end{eqnarray*}
This establishes the next lemma.
 \begin{lemma}\label{lemma1} Let $f:\Omega\to (\tOmega,\rho)$ be harmonic and $\rho$ flat.  Then,  away from the set $\{\mu_f = 0 \}$, 
 \begin{equation}\label{musub}
 \Delta \log |\mu_f| = 0
 \end{equation}
 \end{lemma}
 Next, we set
 \begin{equation}\label{6}
 \sigma^2(z) = \rho^2(f(z)) |f_z(z)|^2 > 0
 \end{equation}
since $\rho>0$ and $f_z\neq 0$.  Also
$\Phi_f = \bar \mu^f \sigma^2$ and so $\log \sigma^2 = \log |\Phi_f|-\log|\mu^f|$ is harmonic away from the isolated set of points where $f_\zbar = 0$. Further, $\log \sigma^2$ is continuous on $\Omega$.  Hence 
 \begin{lemma}\label{lemma2} Let $f:\Omega\to (\tOmega,\rho)$ be harmonic and $\rho$ flat.  Then $\log \sigma$ is harmonic in $\Omega$. \end{lemma}
If $\Omega$ is simply connected (or simply looking locally) we see that $2\log \sigma$ has a harmonic conjugate and so $\sigma^2$ is the modulus of a holomorphic function $\Psi_f$ which does not vanish. Then (\ref{6}) reads as 
 \begin{equation}\label{7}
\bar \mu  \sigma^2(z) = \bar \mu  |\Psi_f| =  \Phi_f  
 \end{equation}
 This gives us a local expression for $|\mu|$
  \begin{lemma}\label{lemma3} Let $f:\Omega\to (\tOmega,\rho)$ be harmonic and $\rho$ flat in a simply connected domain $\Omega$.  Then there is a nonzero holomorphic function $\Psi_f:\Omega\to \IC$ and
 \begin{equation} \label{8}
 |\mu_f| = \frac{|\Phi_f|}{|\Psi_f|}
 \end{equation}
 \end{lemma}
The next two corollaries are obvious consequences.
  \begin{corollary}\label{cor1} Let $f:\Omega\to (\tOmega,\rho)$ be harmonic with $\rho$ flat.  Then $|\mu|$ cannot achieve a local maximum in $\Omega$.  Further,  $|\mu|$ cannot achieve a local minimum unless this value is $0$. \end{corollary}
  \noindent{\bf Remark 1.}   This sort of result cannot be true for more general metrics.  For instance if $\tilde{f}:\Sigma_1\to\Sigma_2$ is a harmonic diffeomorphism between two compact Riemann surfaces, then it must assume its maximal distortion by continuity.  Next $\tilde{f}$ lifts to a diffeomorphism $f:\ID\to\ID$ harmonic with respect to the hyperbolic metric.  This map will, by construction, take its maximum and minimum values of distortion on interior points of the disk (actually any fundamental domain).  Next,  for entire maps we have the following strong rigidity.
    \begin{corollary}\label{cor2} Let $f:\IC \to (\tOmega,\rho)$ be entire,  harmonic with $\rho$ flat.  Then $|\mu| = k <1$ is constant,  $f$ is $\frac{1+k}{1-k}$-quasiconformal and $\tOmega = \IC$. \end{corollary}
\noindent{\bf Proof.}  Evidently the holomorphic function $ \Phi_f/\Psi_f$ is entire and bounded,  and so constant by Liouville's theorem. Since $J(z,f)>0$ at some point,  $|\mu^f|=k<1$ and finally, entire quasiconformal mappings are surjections onto $\IC$.\hfill $\Box$

\medskip

This result should be compared with Wan's result \cite{Wan} which states that an entire harmonic mapping $f:\IH^2\to\IH^2$ of the hyperbolic plane is quasiconformal if and only if the Hopf differential is bounded in the Bloch norm.

\medskip

There is a further consequence of equation (\ref{8}) worth considering.  If $\Omega$ is simply connected,   then both $\Phi_f$ and $\Psi_f$ admit antiderivatives on $\Omega$,  say $\phi$ and $\psi$ respectively.  Set 
\begin{equation}
h(z) = \phi(z) + \overline{\psi}(z)
\end{equation}
Then $h$ is harmonic with respect to the Euclidean metric.  Also $h_z=\Psi_f \neq 0$ and $h_\zbar=\Phi_f$,  so that $K(z,h)=K(z,f)$.  

\medskip

The nonlinear nature of the the tension equation makes it quite difficult to get at exactly what these associated holomorphic function are.  We need this information to effect the isometry to the hyperbolic plane given in Theorem \ref{mainthm}. However,  there is a quite remarkable simplification when we consider the inverse mapping which we now discuss.
  
\section{Harmonic maps and their inverses.}
 
Suppose that $g\in W^{1,2}_{loc}(\tOmega)$ is a quasiconformal  homeomorphism defined on a domain $\tOmega$ in the complex plane with $\Omega=g(\tOmega)$.  We write the Beltrami equation for $g$ as 
\begin{equation}\label{eqn6}
g_\wbar = \mu^g(w) \, g_w
\end{equation}
 Set $f=g^{-1}:\Omega \to (\tOmega,\rho)$ and suppose $f$ is harmonic with respect to $\rho$ (not assumed to be flat yet).  The mapping $f$ is also quasiconformal  and an easy calculation (see \cite[Lemma 10.3.1]{AIM}) shows that $f$ satisfies a Beltrami equation  of the form
\begin{equation}\label{geqn}
f_\zbar = - \mu_g(f(z)) \; \overline{f_z}
\end{equation}
As $f$ is a diffeomorphism, $\mu_g \in C^1(\tOmega)$ and we may differentiate this equation with respect to $z$,  substitute in (\ref{geqn}) to find  (with $\mu_g=\mu$)
\begin{eqnarray*}
f_{\zbar z} & = & - \mu(f(z)) \;  \overline{f_{z \zbar}} - (\mu_w(f(z)) f_z - \mu_\wbar(f(z)) \overline{f_\zbar})\overline{f_z} \\
& = & - \mu(f(z)) \;  \overline{f_{z \zbar}} -  \mu_w(f(z))\,  |f_z|^2 + \mu_\wbar(f(z)) \overline{ \mu(f(z))}  \; |f_z|^2 \\
& = & - \mu(f(z)) \;  \overline{f_{z \zbar}} - \Big( \mu_w(f(z))  - \mu_\wbar(f(z)) \overline{ \mu(f(z))} \Big)  |f_z|^2 
\end{eqnarray*}
Since $f_{z\zbar} + \lambda(f) \, f_z \, f_\zbar = 0$,  $\lambda=(\log \rho)_w$, 
\begin{eqnarray*}
 - \Big( \mu_w(f)  - \mu_\wbar(f) \overline{ \mu(f)} \Big)   |f_z|^2& = &   f_{\zbar z}  + \mu(f) \;  \overline{f_{z \zbar}} \\
 & = &  \lambda(f) \, f_z \, f_\zbar + \mu(f) \;  \overline{\lambda(f) \, f_z \, f_\zbar} \\
  & = & - \lambda(f)  \, \mu(f) \; |f_z|^2 - \mu(g) \;  \overline{\lambda(f)   \, \mu(f) \;  } |f_z|^2 
\end{eqnarray*}
Since $f_z\neq 0$ we obtain the following.
\begin{equation}
 \mu_w(f)  - \mu_\wbar(f) \overline{ \mu(f)}  = - \lambda(f)  \, \mu(f)  - \mu(f) \;  \overline{\lambda(f)   \, \mu(f) \;  } 
\end{equation}
and writing $z=f(w)$ gives us
\begin{equation}\label{eqn9}
 \mu_w  -\overline{ \mu }\; \mu_\wbar    =   \mu\,\Big( \lambda  +   \overline{\mu}\, \overline{\lambda }\Big)
\end{equation}
Of course conversely,  if we have a sufficiently regular  Beltrami  coefficient satisfying this equation we quickly get back to the tension equation for $f$.
We record this in
\begin{theorem}  Let $g:\tOmega \to \Omega$ be quasiconformal,  $\mu = \mu^g$ and suppose $f=g^{-1}:\Omega\to (\tOmega,\rho)$ is harmonic.  Then 
\begin{equation}\label{eqn10}
 \mu_{w}   -\overline{ \mu } \; \mu_{\wbar} =   \mu \,\big( \lambda   +   \overline{\mu }\, \overline{\lambda }\big)
\end{equation}
where $\lambda(w)=(\log \rho)_w(w)$.  Conversely,  if $\mu\in W^{1,1}_{loc}(\tOmega)$,  $\|\mu\|_\infty =k<1$,  satisfies equation (\ref{eqn10}) and if $\mu^g=\mu$,  then $f=g^{-1}:\Omega \to \tOmega$ is harmonic with respect to the metric $\rho(w)|dw|$.
\end{theorem}

\noindent{\bf Example.}  If $\rho(w) = e^{\Im m(w)}>0$,  then $\rho(w)|dw|$ is flat,  $\lambda(w) = i$ and $\mu(w) = \alpha e^{i\Re e(w)/2}$ satisfies equation (\ref{eqn10}).

\medskip

\noindent{\bf Remark 2.}  If $\lambda \equiv 0$,  that is if $\rho$ is constant, then we obtain the class of harmonic mappings with respect to the Euclidean metric and the nonlinear Beltrami equation 
\[ \mu_w  =\overline{ \mu }\; \mu_\wbar \]
 considered in \cite{AIM2}.   
 
\subsection{Flat metrics}

We now assume that $\rho$ determines a flat metric.  We rewrite equation (\ref{eqn9}) and recall $\lambda=(\log \rho)_w$,  so $\bar\lambda=(\log \rho)_\wbar$,  to see
\begin{equation}\label{eqn12}
 \mu_w -  (\log \rho)_w  \; \mu\, = \overline{ \mu }\Big( \mu_\wbar    +   (\log \rho)_\wbar  \; \mu\; \Big)   \end{equation}
Multiply by $e^F$ to get this into   a homogeneous Beltrami equation,   
  \begin{eqnarray*} 
 e^F \mu_w -  e^F (\log \rho)_w  \; \mu\, & = & \overline{ \mu }\Big( e^F \mu_\wbar    +   (\log \rho)_\wbar \; e^F \; \mu\; \Big).
 \end{eqnarray*}
We would like $F_w=-(\log \rho)_w$ and $F_\wbar=(\log \rho)_\wbar$ which would give us
 \begin{eqnarray*} 
 (e^F \mu)_w  & = & \overline{ \mu }\; (e^F \mu)_\wbar    
 \end{eqnarray*}
and this  would make $e^F\mu$ anti-quasiregular.   These conditions on $F$ can be achieved if $\Delta \log\rho = 0$.  We set $F(w)=iv(w)$,  where $v$ is a real valued harmonic conjugate of $\log \rho^{-2}$,  thus $\varphi = -2\log \rho +iv $ is holomorphic.  Then $F=\log\rho+\varphi$ and so $F_\wbar = (\log \rho)_\wbar$ and $F_w=(\log \rho)_z + \varphi' = -(\log \rho)_w$.   We can always construct $v$ on a simply connected domain.  This establishes the next theorem.

\begin{theorem}\label{mainthm}  Let $\Omega$ be simply connected, $f:\Omega \to (\tOmega,\rho)$  harmonic with  $\rho$ defining a  flat metric.  Let $g=f^{-1}:\tOmega\to\Omega$ and let $\mu^g$ be the Beltrami coefficient of $g$. Then,  with $v(w)$ a real valued harmonic conjugate of $\log \rho$ we have  \begin{equation}  \nu(w) = e^{-i v(w)}\overline{\mu^g} (w) \end{equation}
is quasiregular and solves the Beltrami equation
\begin{equation}
\nu_\wbar = \mu(w) \; \nu_w 
\end{equation}
Hence there is a holomorphic $\varphi:\Omega\to \ID$ so that $\nu = \varphi(g)$.
\end{theorem} 
The last part of this theorem is known as Stoilow factorisation,  see \cite{AIM}.  Since quasiregular mappings are open and discrete and since $|\nu(w)|=|\mu^g(w)|$,  this has many obvious consequences .  

\subsection{Harmonic maps with constant distortion}  
Corollary \ref{cor2} shows entire harmonic homeomorphisms with respect to flat metrics have constant distortion.  We take a moment here to investigate the converse. 
Suppose
$ \mu_w -  (\log \rho)_z  \; \mu\, = \overline{ \mu }\Big( \mu_\wbar    +   (\log \rho)_\zbar  \; \mu\; \Big)  $
 and $0\neq |\mu|=k<1$.  Write $\mu(z) = k e^{iv(z)}$,  with $v$ real valued, and calculate that
$ \mu_w    =   i\, \mu \; v_z$, $\mu_\wbar   =    i\,  \mu \; v_\zbar$ so
  \[i\, \mu \; v_z  -  (\log \rho)_z  \; \mu\, = \overline{ \mu }\Big( i\, \mu \; v_\zbar    +   (\log \rho)_\zbar  \; \mu\; \Big)  \]
If $k\neq 0$,  we have $( \log \rho - i\,   v)_z    = - k e^{-iv(z)} (\log \rho + i\,  v )_\zbar $.
But now for a complex function $\overline{h_\zbar} =  \bar h _z$ so either   $(\log \rho - i\,   v)_z   =  (\log \rho + i\,  v )_\zbar = 0$ or $k=1$.  Since we have ruled out the later case,  we must have that $\log \rho + i\,  v $ holomorphic,  $\log \rho$ harmonic and that the metric $\rho(w)|dw|$ flat.
\begin{corollary}  Let $f:\Omega \to (\Omega,\rho)$ be harmonic with constant distortion and not conformal.  Then $\rho$ is flat.
\end{corollary}
Again, the converse is far from true.   There are many harmonic self homeomorphisms of the unit disk with the Euclidean metric which have non-constant distortion.  For instance the Poisson integral extension of almost any homeomorphism $\partial \ID\to\partial \ID$. 
  
\section{Entire harmonic homeomorphisms}

Before establishing Theorem \ref{isothm} we show result Corollary \ref{cor2} cannot be true in complete generality.  

\subsection{Example.}  Let $z=x+iy$ and $\rho(z)=e^{2x}$.  Then $\lambda=(\log \rho)_z = 1$ and $\rho(z)|dz|$ is flat.  We look for a real valued Beltrami coefficient $\mu = \mu(x)$ which solves equation (\ref{eqn10}).  As $\lambda=1$, $\mu_w=\mu_\wbar$ and $\mu=\bar\mu$,  this equation reads as 
$\mu_x  -  \mu   \; \mu_x =   \mu \,\big( 1   +  \mu  \big)$ which we can integrate to find $\mu(1+\mu)^{-2} = a e^x$.  We choose $a=-c$, $c>0$ to get a smooth coefficient
\begin{equation}
\mu = \sqrt{(1+ce^{-x})^2-1}-(1+ce^{-x}), \hskip20pt |\mu|<1
\end{equation}
We can solve the Beltrami equation $f_\zbar = \mu f_z$ as follows.  We write $f(z)=u(x)+iv(y)$.  Then $f_\zbar=u'-v'$ and $f_z=u'+v'$ and we seek $u$ and $v$ so that
\begin{equation}
\frac{u'(x)-v'(y)}{u'(x)+v'(y)} = \sqrt{(1+ce^{-x})^2-1}-(1+ce^{-x})
\end{equation}
Since the right-hand side does not depend on $y$ we may as well assume $v'(y)=1$.  Then
\begin{eqnarray*}
\frac{u'(x)-1}{u'(x)+1} & = & \sqrt{(1+ce^{-x})^2-1}-(1+ce^{-x}) \\
 u'(x)  &=& -1 + \frac{ 2}{2+ce^{-x} - \sqrt{(1+ce^{-x})^2-1}} > 0 
\end{eqnarray*}
If $x<<-1$,  then $u'(x)\sim 1$,  however if $x>>1$,  then $u'(x)\sim \sqrt{\frac{c}{2}} e^{-x/2}$ and so $u(x)$ is bounded and $f$ cannot be surjective.

\subsection{Proof for Theorem \ref{isothm} } 

We first recall the composition formula for the Beltrami coefficients of a quasiconformal mapping.  If $f,g:\IC\to\IC$ are quasiconformal,  then
\begin{equation}
\mu_{f\circ g^{-1}}(g(z)) = \frac{\mu_f(z)-\mu_g(z)}{1-\bar \mu_f(z) \mu_g(z)} \; \frac{\overline{g_z}}{g_z}
\end{equation}
Suppose that $\rho$ is a flat metric density and let $v$ be the harmonic conjugate of $\log \rho$.  If $f$ and $g$ are harmonic with respect to $\rho$,  Theorem \ref{mainthm}  tells us that there are constants $\alpha_f$ and $\alpha_g$ so that $\mu_{f^{-1}}=\alpha_f e^{iv}$ and $\mu_{g^{-1}}=\alpha_g e^{iv}$.  Then the composition formula shows us
\begin{equation}
|\mu_{f^{-1}\circ g}(g^{-1}(z))| = \Big|\frac{(\alpha_f -\alpha_g )e^{iv}}{1-\bar\alpha_f \alpha_g} \; \frac{\overline{g_{z}^{-1}}}{g_{z}^{-1}}\Big| =  \Big|\frac{\alpha_f -\alpha_g}{1-\bar\alpha_f \alpha_g} \Big|
\end{equation}
which provides the isometry with the hyperbolic disk.  Note the surprising fact that $K(z,f^{-1}\circ g)$ is constant under these circumstances.  The equivalence between ${\cal F}_0$ and the disk is then effected by the correspondence
\[ \ID \ni \alpha_f \leftrightarrow \mu_{f^{-1}}= \alpha e^{iv} \leftrightarrow  f \in {\cal F}_0 \]
Continuous dependence on parameters and so forth (see \cite{AIM}) complete the proof. \hfill $\Box$

\bigskip

\noindent G.J. Martin -  Massey University,  Auckland, NZ  \& Magdalen College,  Oxford.

g.j.martin@massey.ac.nz

\end{document}